\documentclass[10pt,doublespacing]{article}
\usepackage{graphicx}
\usepackage{amsmath}
\usepackage{amsxtra}
\usepackage{amsfonts}
\usepackage{amssymb}

\title{A One-Sentence and Truly Elementary Proof of the Basel Problem}
\author{Samuel G. Moreno\thanks{This work was partially supported by  Junta de Andaluc\'{\i}a,
Research Group FQM 0178.}\\{\small{\it Departamento de Matem\'aticas,
Universidad de Ja\'en}}\\{\small{\it 23071 Ja\'en, Spain}}}

\date{}

\begin{document}

\maketitle

\begin{abstract}
By doing a slight change to a beautiful and widely unknown argument by E. L. Stark \cite{St4} we get a candidate to be considered as one of the shortest and most elementary proofs of the celebrated Basel Problem. Furthermore, we give a comprehensive list of references on this topic, displayed in chronological order from Euler to present.
\end{abstract}

\vspace{4.5cm}

The Basel Problem asks for the exact sum  of the reciprocals of the squares of the positive integers. We aim to give a very short and truly elementary proof of this ``golden oldie''.

The first basic idea is to use the well-known trigonometric identity
\begin{eqnarray}\label{e1010}
\frac{1}{2}+\sum_{k=1}^n \cos (k x)=\frac{\sin \big((n+1/2)x\big)}{2 \sin (x/2)},\qquad n\geq 0,
\end{eqnarray}
\noindent which works for all reals $x$ (with the understanding that when $x=2m \pi$, with $m\in \mathbb{Z}$, the right hand side must be taken as its limit value).
The second crucial fact in the proof is the use of the Mean Value Theorem for Integrals, which states that if $f:[a,b]\to\mathbb{R}$ is continuous on  $[a, b]$ and  $g:[a,b]\to\mathbb{R}$
is integrable and nonnegative on  $[a, b]$, then there exists $\xi\in [a, b]$ such that
\begin{eqnarray}
\int_a^b f(x)g(x)\,dx= f(\xi) \int_a^b g(x)\,dx.
\end{eqnarray}

\noindent {\bf Proof.} By first multiplying (\ref{e1010}) by $x^2-2\pi x$, then integrating over $[0,\pi]$ (using integration by parts if necessary)
and finally using the Mean Value Theorem for Integrals, we get

\begin{eqnarray*}
-\frac{\pi^3}{3}+\sum_{k=1}^n \frac{2\pi}{k^2}&=&\int_0^{\pi}\underbrace{(x-2\pi)\frac{x/2}{\sin (x/2)}}_u\underbrace{\vphantom{\frac{x/2}{\sin (x/2)}}\sin \big((n+1/2)x\big)\,dx}_{dv}\\
&=&\underbrace{(x-2\pi)\frac{x/2}{\sin (x/2)}}_u\underbrace{\frac{-\cos\big((n+1/2)x\big)}{n+1/2}}_v\Bigg|_0^{\pi}
-\int_0^{\pi}\underbrace{\frac{-\cos\big((n+1/2)x\big)}{n+1/2}}_v\underbrace{\frac{du}{\vphantom{n+1/2}dx}\,dx}_{du}\\
&=&\frac{-2\pi}{n+1/2}+\frac{\cos\big((n+1/2)\xi_n\big)}{n+1/2}\int_0^{\pi}\frac{du}{dx}\,dx\qquad\qquad\big(\xi_n\in[0,\pi]\big)\\&=&\frac{-2\pi+(u(\pi)-u(0))\cos\big((n+1/2)\xi_n\big)}{n+1/2}=
\frac{-2\pi+(2\pi-\pi^2/2)\cos\big((n+1/2)\xi_n\big)}{n+1/2},
\end{eqnarray*}
\noindent and thus, taking limits when $n$ tends to infinity, we deduce
\begin{eqnarray*}
-\frac{\pi^3}{3}+\sum_{k=1}^{\infty} \frac{2\pi}{k^2}=\lim_{n\to \infty}\frac{-2\pi+(2\pi-\pi^2/2)\cos\big((n+1/2)\xi_n\big)}{n+1/2}=0,
\end{eqnarray*}
\noindent which after some rearrangements gives us
\begin{eqnarray*}
\sum_{k=1}^{\infty} \frac{1}{k^2}=\frac{\pi^2}{6}.
\end{eqnarray*}
{\hfill{$\blacksquare$}}
\noindent

{\bf Remarks}
\begin{enumerate}\leftmargini 1cm

\item Regarding the above use of the Mean Value Theorem for Integrals, let us mention that the function $du/dx$ is nonnegative on $[0,\pi]$ since it is the product of the two increasing functions
 $u_1(x)=x-2\pi$ and $u_2(x)= (x/2)/\sin(x/2)$.

\item Formula (\ref{e1010}) can be easily proved if we first multiply it by $\sin (x/2)$ and then use that $\cos a \sin b= (\sin (a+b)-\sin (a-b))/2$.

\item Our proof is closely connected with the ones in \cite{St4} and \cite{CiNaRuVa}. We came across these two references shortly after a first submission to a journal in December 2014. To tell the truth, our proof is independent of \cite{St4} and \cite{CiNaRuVa}. Moreover, it is the result of a marvelous serendipity when we were preparing for a class on integration methods.

\newpage

\end{enumerate}


\begin{thebibliography}{99}

\bibitem{E1} L. Euler, De Summis Serierum Reciprocarum, {\it Comm. Acad. Sci. Petrop.} {\bf 7} (1734) 123--134.  (Reprinted in Opera Omnia, Ser I.
 14 73--86.)
\bibitem{E2} L. Euler, D\'{e}monstration de la Somme de Cette Suite $1+\frac{1}{4}+\frac{1}{9}+\frac{1}{25}+\frac{1}{36}+{\rm etc}$,
{\it Journ. Lit. d'Allemagne, de Suisse et du Nord} {\bf 2} (1743) 115--127.  (Reprinted in Opera Omnia, Ser I. 14 177--186.)
\bibitem{E3} L. Euler, De Summatione Serierum in hac Forma Contentarum $a/1+a^2/4+a^3/9+a^4/16+a^5/25+a^6/36+{\rm etc}$,
{\it Memories de l'academie des Sciences de St. Petersbourg} {\bf 3} (1811) 26--42. (Reprinted in Opera Omnia, Ser I. 16 117--138.)
\bibitem{Stielt} T. J. Stieltjes, Table des Valeurs des Sommes $S_k=\sum_1^\infty n^{-k}$, {\it Acta Math.} {\bf 10} (1987) 299--302.
\bibitem{Stack} P. St\"{a}ckel, Eine Vergessene Abhandlung Leonhard Eulers \"{u}ber die Summe der Reziproken Quadrate der Nat\"{u}rlichen Zahlen, {\it Bibliotheca Mathematic 3} {\bf 8} (1907) 37--60.
 (Available at http://archiv.ub.uni--heidelberg.de/volltextserver/13424/1/staeckel{\_ }Euler{\_ }Abh.pdf.)
\bibitem{KnoSch} K. Knopp, I. Schur, \"{U}ber die Herleitung der Gleichung $\sum_{n=1}^{\infty}1/n^2=\pi^2/6$, {\it Arch. der Math. u. Phys.(3)} {\bf 27} (1918) 174--176.
\bibitem{Titch} E. C. Titchmarsh, A Series Inversion Formula, {\it Proc. London Math. Soc. (2)} {\bf 26} (1927) 1--11.
\bibitem{Est} T. Estermann, Elementary Evaluation of $\zeta(2k)$, {\it J. London Math. Soc.} {\bf 22} (1947) 10--13.
\bibitem{Kuo} H. T. Kuo, A Recurrence Formula for $\zeta(2n)$, {\it Bull. Amer. Math. Soc.} {\bf 55} (1949) 573--574.
\bibitem{Will1} G. T. Williams, A New Method of Evaluating $\zeta(2n)$, {\it Amer. Math. Monthly} {\bf 60} (1953) 19--25.
\bibitem{Mat} Y. Matsuoka, An Elementary Proof of the Formula $\sum_{k=1}^\infty1/k^2=\pi^2/6$, {\it Amer. Math. Monthly} {\bf 68} (1961) 485--487.
\bibitem{YY} A. M. Yaglom, I. M. Yaglom, {\it Challenging Mathematical Problems with Elementary Solutions}, Vol. II (Problem 145), Holden-Day, San Francisco, 1967.
\bibitem{St1} E. L. Stark, Another Proof of the Formula $\sum 1/k^2=\pi^2/6$, {\it Amer. Math. Monthly} {\bf 76} (1969) 552--553.
\bibitem{Hol} F. Holme, En Enkel Beregning av $\sum_{k=1}^\infty 1/k^2$, {\it Nordisk Mat. Tidskr.} {\bf 18} (1970) 91--92.
\bibitem{St2} E. L. Stark, $1-\frac{1}{4}+\frac{1}{9}-\frac{1}{16}+\cdots=\frac{\pi^2}{12}$, {\it Praxis Math.} {\bf 12}  (1970) 1--3.
\bibitem{SkSel} I. S. Skau, E. S. Selmer, Noen Anvendelser av Finn Holmes Methode for Beregning av  $\sum_{k=1}^\infty \frac{1}{k^2}$, {\it Nordisk Mat. Tidskr.} {\bf 19} (1971) 120--124.
\bibitem{Will2} K. S. Williams, On $\sum_{n=1}^\infty (1/n^{2k})$, {\it Math. Mag.} {\bf 44} (1971) 273--276.
\bibitem{Gies} D. P. Giesy, Still Another Elementary Proof that $\sum 1/k^2=\pi^2/6$, {\it Math. Mag.} {\bf 45} (1972) 148--149.
\bibitem{P} I. Papadimitriou, A Simple Proof of the Formula $\sum_{k=1}^\infty k^{-2}=\pi^2/6$, {\it Amer. Math. Monthly} {\bf 80} (1973) 424--425.
\bibitem{Ap1} T. M. Apostol, Another Elementary Proof of Euler's Formula for $\zeta(2n)$, {\it Amer. Math. Monthly} {\bf 80} (1973) 425--431.
\bibitem{Ayo} R. Ayoub, Euler and the Zeta Function, {\it Amer. Math. Monthly} {\bf 81} (1974) 1067--1086.
\bibitem{St3} E. L. Stark, The Series $\sum_{k=1}^\infty k^{-s}$, $s=2,\,3,\, 4,\ldots$ Once More, {\it Math. Mag.} {\bf 47} (1974) 197--202.
\bibitem{Bern} B. C. Berndt, Elementary evaluation of  $\zeta(2n)$, {\it Math. Mag.} {\bf 48} (1975) 148--154.
\bibitem{Chen} M. P. Chen, An Elementary Evaluation of $\zeta(2n)$, {\it Chinese J. Math.} {\bf 3} (1975) 11--15.
\bibitem{St4} E. L. Stark, Application of a Mean Value Theorem for Integrals to Series Summation, {\it Amer. Math. Monthly} {\bf 85} (1978) 481--483.
\bibitem{Byu} B. I. Kim, Application of a Mean Value Theorem for Integrals to Prove the Formula $\displaystyle{\sum \frac{1}{K^2}=\frac{\pi^2}{6}}$, {\it J. Korea Soc. Math. Educ.} {\bf 17} (1979) 33--35. (This seems to be a clone of \cite{St4}.)
\bibitem{Beuk} F. Beukers, A Note on the Irrationality of $\zeta(2)$ and $\zeta(3)$, {\it Bull. London Math. Soc.} {\bf 11} (1979) 268--272.
\bibitem{Ap2} T. M. Apostol, A Proof that Euler Missed: Evaluating $\zeta(2)$ the Easy Way, {\it Math. Intelligencer} {\bf 5} (1983) 59--60.
\bibitem{ErDud} P. Erd\"{o}s, U. Dudley,  Some Remarks and Problems in Number Theory Related to the Work of Euler, {\it Math. Mag.} {\bf 56} (1983) 292--298.
\bibitem{Klin} M. Kline, Euler and Infinite Series, {\it Math. Mag.} {\bf 56} (1983) 307--314.
\bibitem{Kim} G. Kimble, Euler's Other Proof, {\it Math. Mag.} {\bf 60} (1987) 282.
\bibitem{Ch} B. R. Choe, An Elementary Proof of $\sum_{n=1}^\infty 1/n^2=\pi^2/6$, {\it Amer. Math. Monthly} {\bf 94} (1987) 662--663.
\bibitem{LeSchu} B. Leonard, H. S. Shultz, A Computer Verification of a Pretty Mathematical Result, {\it Math. Gazette} {\bf 72} (1988) 7--10.
\bibitem{Sh} N. Shea, Summing the series $\frac{1}{1^2}+\frac{1}{2^2}+\frac{1}{3^2}+\ldots$, {\it Math. Spectrum} {\bf 21} (1989) 49--55.
\bibitem{Russ} D. C. Russell, Another Eulerian-Type Proof, {\it Math. Mag.} {\bf 64} (1991) 349.
\bibitem{Kal} D. Kalman, Six Ways to Sum a Series, {\it College Math. J.} {\bf 24} (1993) 402--421.
\bibitem{BeuKolCal} F. Beukers, J. A. C. Kolk, E. Calabi, Sums of Generalized Harmonic Series and Volumes, {\it Nieuw Arch. Wisk.} {\bf 11} (1993)  217--224.
\bibitem{Kor} R. A. Kortram, Simple Proofs for $\sum\limits^\infty_{k=1} \frac{1}{k^2} = \frac{\pi^2} {6}$ and
$\sin x = x\prod\limits^\infty_{k=1} \bigg(1 - \frac{x^2}{k^2\pi^2} \bigg)$, {\it Math. Mag.} {\bf 69} (1996) 122--125.
\bibitem{MckinTuck} M. McKinzie, C. Tuckey, Hidden Lemmas in Euler's Summation of the Reciprocals of the Squares, {\it Arch. Hist. Exact Sci.} {\bf 51}
 (1997) 29--57.
\bibitem{AigZieg} M. Aigner, G. M. Ziegler, {\it Proofs from THE BOOK} (Chapter 7), third edition, Springer, Berlin, 1998.
\bibitem{Robb} N. Robbins, Revisiting and Old Favorite: $\zeta (2m)$, {\it Math. Mag.} {\bf 72} (1999) 317--319.
\bibitem{Chap} R. Chapman, Evaluating $\zeta (2)$, (1999). (Available at

http://www.uam.es/personal{\_}pdi/ciencias/cillerue/Curso/zeta2.pdf.)
\bibitem{Balan} E. P. Balanzario, M\'etodo Elemental para la Evaluaci\'on de la Funci\'on Zeta de Riemann en los Enteros Pares, {\it Miscel\'anea Mat.}
 {\bf 33} (2001) 31--41.
\bibitem{Sang} C. J. Sangwin, An Infinite Series of Surprises, (2001). (Available at http://venus.unive.it/pasinigi/teaching/Baselprob.pdf.)
\bibitem{Huy} D. Huylebrouck, Similarities in Irrationality Proofs for $\pi$, $\ln 2$, $\zeta(2)$, and $\zeta(3)$, {\it Amer. Math. Monthly} {\bf 108}
 (2001) 222--231
\bibitem{Hof} J. Hofbauer, A Simple Proof of  $1+\frac{1}{2^2}+\frac{1}{3^2}+\cdots=\frac{\pi^2}{6}$ and Related Identities, {\it Amer. Math. Monthly}
 {\bf 109} (2002) 196--200.
\bibitem{Lord} N. Lord, Yet Another Proof that $\sum \frac{1}{n^2}=\frac{\pi^2}{6}$, {\it Math. Gazette} {\bf 86} (2002) 477--479.
\bibitem{Har} J. D. Harper, Another Simple Proof of  $1+\frac{1}{2^2}+\frac{1}{3^2}+\cdots=\frac{\pi^2}{6}$, {\it Amer. Math. Monthly} {\bf 110} (2003) 540--541.
\bibitem{Stopp} J. Stopple, A Primer of Analytic Number Theory. From Pythagoras to Riemann, Chapter 6, {\it Cambridge University Press, New York}, 2003.
\bibitem{Sand1} C. E. Sandifer, Estimating the Basel Problem, (2003). (Available at http://eulerarchive.maa.org/hedi/HEDI-2003-12.pdf.)
\bibitem{Os} T. J. Osler, Finding $\zeta (2p)$ from a Product of Sines, {\it Amer. Math. Monthly} {\bf 111} (2004) 52--54.
\bibitem{Kok} K. P. Kokhas, Sum of Inverse Squares, {\it Mathematicheskoe Prosveshenie} {\bf 8} (2004) 142--163.
\bibitem{Tsu} H. Tsumura, An Elementary Proof of  Euler's Formula for $\zeta(2m)$, {\it Amer. Math. Monthly} {\bf 111} (2004) 430--431.
\bibitem{Sand2} C. E. Sandifer, Basel Problem with Integrals, (2003). (http://eulerarchive.maa.org/hedi/HEDI-2004-03.pdf.)
\bibitem{BourFuYor} P. Bourgade, T. Fujita, M. Yor, Euler's Formulae for $\zeta (2n)$ and Products of Cauchy Variables, {\it Electron. Commun. Probab.}
{\bf 12} (2007) 73--80.
\bibitem{BraDanSan} Euler at 300. An Appreciation. Edited by R. E. Bradley, L. A. D'Antonio and C. E. Sandifer,
{\it Mathematical Association of America, Washington, DC}, 2007.
\bibitem{Sand3} C. E. Sandifer, Euler's Greatest Hits, (2007). (http://eulerarchive.maa.org/hedi/HEDI-2007-02.pdf.)
\bibitem {Bred} M. Brede, Eulers Identit\"{a}ten f\"{u}r die Werte von $\zeta (2n)$, {\it Mathematische Semesterberichte} {\bf 54} (2007) 135--140.
\bibitem{Iv} M. Ivan, A Simple Solution to Basel Problem, {\it Gen. Math.} {\bf 16} (2008) 111--113.
\bibitem{P} M. Passare, How to Compute $\sum 1/n^2$ by Solving Triangles, {\it Amer. Math. Monthly} {\bf 115} (2008) 745--752.
\bibitem{EvRotWa} G. Everest, C. R\"{o}ttger, T. Ward, The Continuing Story of Zeta, {\it Math. Intelligencer} {\bf 31} (2009) 13--17.
\bibitem{Sila} Z. K. Silagadze, Sums of Generalized Harmonic Series for Kids from Five to Fifteen, (2010). (Available at http://arxiv.org/pdf/1003.3602v1.pdf.)
\bibitem{Mar} T. Marshall, A Short Proof of $\zeta(2)=\pi^2/6$, {\it Amer. Math. Monthly} {\bf 117} (2010) 352--353.
\bibitem{Wast} J. W\"{a}stlund, Summing Inverse Squares by Euclidean Geometry, (2010). (Available at http://www.math.chalmers.se/$\sim$wastlund/Cosmic.pdf.)
\bibitem{Lev} P. Levrie, Lost and Found: An Unpublished  $\zeta(2)$-Proof, {\it Math. Intelligencer} {\bf 33} (2011) 29--32.
\bibitem{Hich} M. D. Hirschhorn, A Simple Proof that $\zeta(2)=\frac{\pi^2}{6}$, {\it Math. Intelligencer} {\bf 33} (2011) 81--82.
\bibitem{AmCaFer} E. De Amo, M. D\' {i}az Carrillo, J. Fern\'andez--S\'anchez, Another Proof of Euler's Formula for $\zeta (2k)$, {\it Proc. Amer. Math. Soc.}
 {\bf 139} (2011) 1441--1444.
\bibitem{Pac} L. Pace, Probabilistically Proving  that $\zeta(2)=\pi^2/6$, {\it Amer. Math. Monthly} {\bf 118} (2011) 641--643.
\bibitem{KalMckin} D. Kalman, M. Mckinzie, Another Way to Sum a Series: Generating Functions, Euler, and the Dilog Function, {\it Amer. Math. Monthly}
{\bf 119} (2012) 42--51.
\bibitem{Ben} D. Benko, The Basel Problem as a Telescoping Series, {\it College Math. J.} {\bf 43} (2012) 244--250.
\bibitem{Dan} D. Daners, A Short Elementary Proof of $\sum 1/k^2=\pi^2/6$, {\it Math. Mag.} {\bf 85} (2012) 361--364.
\bibitem{Li1} F. M. S. Lima, Another Elementary Proof of $\sum_{n\geq 1} 1/n^2=\pi^2/6$ and a Recurrence Formula for  $\zeta(2k)$, (2012).
(Available at http://arxiv.org/pdf/1109.4605.pdf.)
\bibitem{CiNaRuVa} O. Ciaurri, L. M. Navas, F. J. Ruiz, J. L. Varona, A Simple Computation of $\zeta(2k)$ by using Bernoulli Polynomials and a Telescoping Series, (2012). (Available at http://arxiv.org/pdf/1209.5030v1.pdf.)
\bibitem{Dalai} M. Dalai, How Would Riemann Evaluate $\zeta (2n)$, {\it Amer. Math. Monthly} {\bf 120} (2013) 169--171.
\bibitem{Ritell} D. Ritelli, Another Proof of $\zeta(2)=\frac{\pi^2}{6}$ Using Double Integrals, {\it Amer. Math. Monthly} {\bf 120} (2013) 642--645.
\bibitem{Pat} S. Patr\`{i}, Sum of the Generalized Harmonic Series with Even Natural Exponents, {\it Rend. Mat. Appl. (7)} {\bf 33} (2013) 19--26.
\bibitem{Jame} T. Jameson, Another Proof that $\zeta(2)=\pi^2/6$ Via Double Integration, {\it Math. Gazette} {\bf 97} (2013) note 97.44.
\bibitem{JameLord} G. J. O. Jameson, N. Lord, Evaluation of $\sum_{n=1}^{\infty}\frac{1}{n^2}$ by a Double Integral, {\it Math. Gazette} {\bf 97} (2013).
\bibitem{Gay} J. Gayo, O Problema que Tornou Euler Famoso,  Disserta\c{c}\~{a}o de Mestrado, (2013).
(Available at http://repositorio.utfpr.edu.br/jspui/handle/1/579?mode=full.)
\bibitem{XuZho} H. Xu, J. Zhou, The Connection between the Basel Problem and a Special Integral, (2013). (Available at   http://arxiv.org/abs/1307.8278.)
\bibitem{Sull} B. W. Sullivan, The Basel Problem. Numerous results, (2013). (Available at http://math.cmu.edu/$\sim$bwsulliv/basel-problem.pdf.)
\bibitem{Kraus} R. M. Krause, $\zeta(2)$ Once Again, {\it Amer. Math. Monthly} {\bf 121} (2014) 353--354.
\bibitem{XuZh} H. Xu, J. Zhou, The Connection between the Basel Problem and a Special Integral, {\it Applied Mathematics} {\bf 5} (2014) 2570--2584.
\bibitem{Sev} H. Lundmark et al., Different Methods to Compute $\sum_{n=1}^{\infty}\frac{1}{n^2}$, (2010-2014).


(http://math.stackexchange.com/questions/8337/different-methods-to-compute-sum-limits-n-1-infty-frac1n2.)
\bibitem{Swan} J. W. H. Swanepoel, On a Generalization of a Theorem by Euler, {\it J. Number Theory} {\bf 149} (2015) 46--56.

\end{thebibliography}
\end{document}